\documentclass{article}
\usepackage{amsmath}
\usepackage{amssymb}
\usepackage{amsthm}
\usepackage{array}
\usepackage{booktabs}
\usepackage{float}
\usepackage[colorlinks]{hyperref}
\usepackage{verbatim}
\providecommand{\keywords}[1]{\textbf{\textit{Keywords:}} #1}

\title{A symmetric chain decomposition of $L(6,n)$}
\author{Xiangdong Wen\\
Wolfram Research}
\date{}

\newtheorem{theorem}{Theorem}

\begin{document}
\maketitle

\begin{abstract}
Stanley conjectured that Young's lattice $L(m,n)$ has a symmetric chain decomposition. In this paper we give an explicit symmetric chain decomposition of $L(6,n)$ for all $n$. The decomposition is encoded by parametric chain families in Wolfram Language. The proof has two parts. First, the chain-family program produces disjoint saturated symmetric chains. Second, the total weight of all chains is shown to equal the generating function of $L(6,n)$. The weight identity is reduced to an exact rational-function identity with $2324$ grouped terms, then verified by a recursive pole-vanishing argument and a finite exact grid check.
\end{abstract}

\keywords{Young's lattice, poset, symmetric chain decomposition, computer proof}

\section{Introduction}

For positive integers $m$ and $n$, Young's lattice $L(m,n)$ is the partially ordered set of $m$-tuples
\[
(a_1,a_2,\ldots,a_m), \qquad 0 \le a_1 \le a_2 \le \cdots \le a_m \le n,
\]
ordered coordinatewise. The rank of $(a_1,\ldots,a_m)$ is $\sum_{r=1}^m a_r$. A chain
\[
\mathbf{v}_1 < \mathbf{v}_2 < \cdots < \mathbf{v}_s
\]
is saturated if it skips no rank, and symmetric if
\[
\operatorname{rank}(\mathbf{v}_1)+\operatorname{rank}(\mathbf{v}_s)=mn.
\]
A symmetric chain decomposition (SCD) of a graded poset is a partition into saturated symmetric chains.

Stanley \cite{Stanley} conjectured that $L(m,n)$ has an SCD for all $m,n$. Explicit constructions are known for small $m$; see \cite{Lind,West,Wen04,Kathy,Doron,Wen24}. Related structural work on Young's lattice and nearby posets includes Dhand's tropical decomposition approach \cite{Dhand} and Zhong's construction for weak-composition rank sequences \cite{Zhong}. The purpose of this paper is to record a complete explicit construction for $L(6,n)$.

The main result is the following.

\begin{theorem}
For every nonnegative integer $n$, the chains defined by the chain-family program form a symmetric chain decomposition of $L(6,n)$.
\end{theorem}

The proof is computer-assisted but exact. The decomposition itself is given by explicit formulas in the Wolfram Language source linked above, and the covering identity is proved by exact symbolic calculations described in Section \ref{sec:proof}.

\section{Chain Families}

The chain-family program organizes the decomposition into $91$ parametric chain families. Their names are grouped by prefix, and the groups are the natural units of the construction.

There are two levels of notation in the first block of the construction.
\begin{itemize}
\item \texttt{p01}, \dots, \texttt{p11} denote parallel chains, that is, parametric saturated chains arranged into rectangular families.
\item \texttt{c01}, \dots, \texttt{c11} denote the symmetric chains obtained from those parallel chains by applying \texttt{symmetricChain}.
\end{itemize}
All remaining named families, namely the \texttt{cc}, \texttt{cd}, \texttt{cf}, \texttt{ce}, \texttt{cu}, \texttt{cx}, \texttt{cz}, \texttt{cs}, and \texttt{cr} blocks, are already given directly as symmetric chain families in the program.

The families are grouped by prefix in Table \ref{tab:family-groups}.

\begin{table}[H]
\centering
\begin{tabular}{lc}
\toprule
Group & Number of families \\
\midrule
\texttt{p} & 11 \\
\texttt{cc} & 3 \\
\texttt{cd} & 3 \\
\texttt{cf} & 13 \\
\texttt{ce} & 17 \\
\texttt{cu} & 2 \\
\texttt{cx} & 3 \\
\texttt{cz} & 4 \\
\texttt{cs} & 20 \\
\texttt{cr} & 15 \\
\midrule
Total & 91 \\
\bottomrule
\end{tabular}
\caption{Chain-family groups in the Wolfram Language construction.}
\label{tab:family-groups}
\end{table}

The utilities
\[
\texttt{sym2DChain},\quad \texttt{symmetricChain},\quad \texttt{L2}
\]
play the same role as in the $L(5,n)$ construction: they convert rectangular or two-parameter collections of parallel chains into saturated symmetric chains by recursively taking perimeters. Here two chains are called parallel when, with some parameters fixed and one or two parameters varying, they form a rectangular array of saturated chains with the same combinatorial shape.
This is exactly how \texttt{c01}, \dots, \texttt{c11} are produced from \texttt{p01}, \dots, \texttt{p11}. The other families are written directly as symmetric chains, with more parameters and a larger range of affine constraints on those parameters.

A representative example is the pair \texttt{p01}/\texttt{c01}. The function \texttt{p01} is a saturated chain depending on the parameters $n,i,j,k,h,p,q$. The corresponding family \texttt{c01} is obtained by applying \texttt{symmetricChain} to the rectangular array of chains obtained when the parameter pair $(p,q)$ ranges over the two-dimensional set \texttt{L2[k,k]}. Thus \texttt{p01} is the parallel input and \texttt{c01} is the resulting family of symmetric chains. For instance,
\[
\texttt{p01[2,0,0,0,0,0,0]}
\]
is the chain
\[
\begin{aligned}
&(0,0,0,0,1,1)<(0,0,0,1,1,1)<(0,0,1,1,1,1)\\
&<(0,1,1,1,1,1)<(1,1,1,1,1,1)<(1,1,1,1,1,2)\\
&<(1,1,1,1,2,2)<(1,1,1,2,2,2)<(1,1,2,2,2,2),
\end{aligned}
\]
which is already a saturated chain in $L(6,2)$ and sits inside the larger rectangular family from which \texttt{c01} is assembled.

A representative direct symmetric family is \texttt{cc0}. Here the program defines \texttt{cc0[n,h,i,j,k]} directly as a concatenation of tables describing a single symmetric chain, and then forms the full family by ranging over all admissible parameter values and flattening the resulting list. For example,
\[
\texttt{cc0[1,0,0,0,0]}
\]
is the symmetric chain
\[
\begin{aligned}
&(0,0,0,0,0,0)<(0,0,0,0,0,1)<(0,0,0,0,1,1)<(0,0,0,1,1,1)\\
&<(0,0,1,1,1,1)<(0,1,1,1,1,1)<(1,1,1,1,1,1),
\end{aligned}
\]
whose endpoint ranks sum to $0+6=6=6\cdot 1$.
Thus \texttt{cc0} has no intermediate parallel-chain stage; it is already part of the final symmetric decomposition.

The formulas are explicit, but printing all $91$ families in the main body would obscure the construction. The exact chain formulas are supplied in the linked Wolfram Cloud source. The important points for the proof are:
\begin{itemize}
\item each family consists of vectors in $L(6,n)$,
\item each \texttt{p}-family is a family of parallel saturated chains,
\item each \texttt{c}-family is obtained from a corresponding \texttt{p}-family by the helper routine \texttt{symmetricChain},
\item each remaining family is already presented as a symmetric chain family.
\end{itemize}

\section{From Chains to Weights}
\label{sec:weights}

For
\[
\mathbf{a}=(a_1,a_2,a_3,a_4,a_5,a_6)\in L(6,n),
\]
define the weight
\[
w(\mathbf{a})=x_0^{n-a_6}x_1^{a_6-a_5}x_2^{a_5-a_4}x_3^{a_4-a_3}x_4^{a_3-a_2}x_5^{a_2-a_1}x_6^{a_1}.
\]
This is exactly the weight definition used in the Wolfram Language verification.

The total weight of $L(6,n)$ over all $n$ is
\[
\sum_{n\ge 0}\ \sum_{\mathbf{a}\in L(6,n)} w(\mathbf{a})
= \frac{1}{(1-x_0)(1-x_1)(1-x_2)(1-x_3)(1-x_4)(1-x_5)(1-x_6)}.
\]
So it is enough to prove that the total weight of all chain families equals this generating function.

The point of this weight map is that distinct vectors of $L(6,n)$ give distinct monomials. Hence equality of generating functions is exactly the statement that every element of $L(6,n)$ occurs with coefficient $1$ in the union of the chain families, that is, the chains cover $L(6,n)$ without overlap or omission.

The combined weight/check program carries out this conversion in three stages.
\begin{enumerate}
\item The routine \texttt{chainWeight} converts a chain family into a list of weighted monomials with parameter conditions.
\item The routine \texttt{sp} splits the resulting conditions into disjoint regions by resolving the \texttt{Min} expressions introduced by the parametrizations.
\item The symbols \texttt{tw01}, \texttt{tw02}, \dots, \texttt{twr25} collect all weighted contributions.
\end{enumerate}

Then the program forms
\[
\texttt{ss = Names["tw*"]},
\]
extracts all split pieces, and applies \texttt{NonNegativeLatticeSum} to each piece. This yields $2324$ grouped rational terms. Therefore the covering identity reduces to
\[
\sum_{r=1}^{2324} T_r
=
\frac{1}{(1-x_0)(1-x_1)(1-x_2)(1-x_3)(1-x_4)(1-x_5)(1-x_6)}.
\]
Equivalently,
\[
\sum_{r=1}^{2324} T_r-\frac{1}{(1-x_0)(1-x_1)(1-x_2)(1-x_3)(1-x_4)(1-x_5)(1-x_6)}=0.
\]
Since the chain families are already saturated and symmetric by construction, this identity is the only remaining step needed to prove that they form an SCD.

\section{Proof of the Weight Identity}
\label{sec:proof}

The difference above is exact but large, so we do not prove it by a single global \texttt{Together}. Instead we use an exact two-stage argument on the expanded atomic terms.

\subsection*{Stage 1: Vanishing of all denominator factors}

Each of the $2324$ grouped terms is expanded into atomic rational terms. After adding the correction term
\[
-\frac{1}{(1-x_0)(1-x_1)(1-x_2)(1-x_3)(1-x_4)(1-x_5)(1-x_6)},
\]
this produces $6800$ atomic terms.

We then consider every distinct irreducible denominator factor appearing in those terms. There are $532$ such factors in the expanded expression. For each factor $f$, we prove recursively that the total residue of the expression along $f=0$ vanishes. The Wolfram Language proof returns
\[
\mathrm{TotalFactors}=532,\qquad
\mathrm{ClosedCount}=532,\qquad
\mathrm{OpenCount}=0.
\]
Hence every denominator factor has zero pole, so the total expression is a polynomial. This is the exact point where the rational identity becomes a polynomial identity.

\subsection*{Stage 2: Exact polynomial zero test}

Once polynomiality is known, the Wolfram Language proof gives exact degree bounds in each variable:
\[
\begin{aligned}
\deg_{x_0}&\le 1, & \deg_{x_1}&\le 2, & \deg_{x_2}&\le 3,\\
\deg_{x_3}&\le 2, & \deg_{x_4}&\le 2, & \deg_{x_5}&\le 1,\\
\deg_{x_6}&\le 2. &&
\end{aligned}
\]

These bounds imply that it is enough to check the polynomial on a Cartesian grid of size
\[
(1+1)(2+1)(3+1)(2+1)(2+1)(1+1)(2+1)=1296.
\]
The final exact check evaluates the expression at the grid
\[
\begin{aligned}
x_0&\in\{1297,1298\},\\
x_1&\in\{1299,1300,1301\},\\
x_2&\in\{1302,1303,1304,1305\},\\
x_3&\in\{1306,1307,1308\},\\
x_4&\in\{1309,1310,1311\},\\
x_5&\in\{1312,1313\},\\
x_6&\in\{1314,1315,1316\},
\end{aligned}
\]
and finds zero at every point. In particular,
\[
\mathrm{GridSize}=1296,\qquad \mathrm{BadPointCount}=0.
\]

Therefore the polynomial is identically zero. So the rational difference vanishes identically, and the total weight of the chain families is exactly
\[
\frac{1}{(1-x_0)(1-x_1)(1-x_2)(1-x_3)(1-x_4)(1-x_5)(1-x_6)}.
\]
Because the weight map is injective on lattice elements, this means that every element of $L(6,n)$ appears exactly once in the union of the chain families. Since those chains are already saturated and symmetric by construction, they form a symmetric chain decomposition of $L(6,n)$, proving the theorem.

\section*{Code Availability}

The construction and verification are recorded in Wolfram Cloud objects:
\begin{itemize}
\item chain-family program:
\href{https://www.wolframcloud.com/obj/ad584764-8bbb-44b9-8736-53bbd6da2703}{Wolfram Cloud source 1}
\item combined weight/check program:
\href{https://www.wolframcloud.com/obj/106e2b9d-1ad8-4a79-a853-dbc265445505}{Wolfram Cloud source 2}
\item proof programs:
\href{https://www.wolframcloud.com/obj/5406c3c6-00dd-481c-9dce-7f2f6988feb0}{Wolfram Cloud source 3},
\href{https://www.wolframcloud.com/obj/8bf07f29-11a3-480b-8979-3cd4a2227740}{Wolfram Cloud source 4},
\href{https://www.wolframcloud.com/obj/4228238b-f011-491b-8c09-96a78d8e36d9}{Wolfram Cloud source 5}
\item the $2324$ split conditions:
\href{https://www.wolframcloud.com/obj/f86200a5-8c8d-4cab-9b21-b67fc07f60a7}{Wolfram Cloud source 6}
\item the $2324$ lattice-sum outputs:
\href{https://www.wolframcloud.com/obj/54923749-02c6-4deb-9a10-654e2902b248}{Wolfram Cloud source 7}
\end{itemize}

\section*{Acknowledgment}

I am indebted to Dr. Kathy O'Hara for discussions on symmetric chain decompositions and to Doron Zeilberger for encouragement to write up these constructions. I am also grateful to collaborators and colleagues who discussed the computational side of the project.

\end{document}